\documentclass[a4paper,10pt,draft]{article}
\usepackage{amssymb, amsmath, amsthm, latexsym, verbatim}

\newcommand \G{\Gamma}
\newcommand \n{\nabla}
\newcommand \om{\omega}
\newcommand \Om{\Omega}
\newcommand \M{\mathcal M}
\newcommand \la{\lambda}
\newcommand \id{\mathrm{id}}
\newcommand \Rn{\mathbb R^n}
\newcommand \Ro{\mathbb R^8}

\newcommand \RS{\mathbb R^{7}}
\newcommand \RR{\mathbf R}

\newcommand \rk{\mathrm{rk}}
\newcommand \Ker{\mathrm{Ker}}
\newcommand \Ric{\mathrm{Ric}}

\newcommand \Cl{\mathrm{Cl}}
\newcommand \Span{\mathrm{Span}}
\newcommand \Range{\mathrm{Range}}
\newcommand \Tr{\mathrm{Tr}}
\newcommand \Hom{\mathrm{Hom}}

\newcommand \A{\mathfrak{A}_0}
\newcommand \K{\mathfrak{K}_0}
\newcommand \W{\|X \wedge Y\|^2}
\renewcommand \a{\alpha}
\renewcommand \b{\beta}

\newcommand \<{\langle}
\renewcommand \>{\rangle}

\theoremstyle{plain}

\newtheorem*{theorem*}{Theorem}
\newtheorem*{corollary*}{Corollary}
\newtheorem{lemma}{Lemma}

\newtheorem*{proposition*}{Proposition}
\newtheorem{fact}{Fact}
\newtheorem*{namedtheorem}{\theoremname}
\newcommand{\theoremname}{te}

\theoremstyle{definition}

\newtheorem*{definition*}{Definition}

\theoremstyle{remark}
\newtheorem*{remark}{Remark}
\newtheorem{example}{Example}

\begin{document}

\title{Riemannian manifolds of dimension $7$ whose skew-symmetric curvature operator has constant
eigenvalues}

\author{Y.Nikolayevsky\thanks{Work supported by the ARC Discovery grant DP0342758}}

\date{}

\maketitle

\begin{abstract} A Riemannian manifold $M^n$ is called IP, if, at every point $x \in M^n$,
the eigenvalues of its skew-symmetric curvature operator $R(X, Y)$ are the same,
for every pair of orthonormal vectors $X, Y \in T_xM^n$. In \cite{GLS, G2, IP} it was
shown that for all $n \ge 4$, except $n=7$, any IP manifold either has constant curvature,
or is a warped product, with some specific function, of a line and a space of constant
curvature. We extend this result to the case $n = 7$, and also discuss $3$-dimensional
IP manifolds.
\end{abstract}

\section{Introduction}\label{s:intro}

An {\it algebraic curvature tensor} $R$ in a Euclidean space $\Rn$ is a $(3, 1)$
tensor having the same symmetries as the curvature tensor of a Riemannian manifold.
Given an algebraic curvature tensor $R$, there defined a quadrilinear functional on $\Rn$
by $R(X,Y,Z,W) = \<R(X,Y)W, Z\>$.
For any pair of vectors $X, Y \in \Rn, \, R(X, Y)$ is a skew-symmetric endomorphism
of $\Rn$. One has $R(Y, X) = - R(X, Y)$, and, in particular, $R(X, Y) =0$ when $X \parallel Y$.
For any oriented two-plane $\pi \in G^+(2,n)$, there is a well-defined endomorphism
$R(\pi)$ of $\Rn, \quad R(\pi) = \|X \wedge Y\|^{-1} R(X, Y)$, where $(X, Y)$ is any
oriented pair of vectors spanning $\pi$, and
$\|X \wedge Y\| = (\|X\|^2\|Y\|^2 - \<X, Y\>^2)^{1/2}$.

\begin{definition*} An algebraic curvature tensor $R$ is called
{\it IP}, if the eigenvalues of $R(\pi)$ are the same for all $\pi \in G^+(2,n)$.
A Riemannian manifold $M^n$ is called {\it IP}, if its curvature tensor at every point
is IP (the eigenvalues may depend on a point).
\end{definition*}

For an IP algebraic curvature tensor $R$, its \textit{rank} is the rank of any
of the $R(\pi)$'s.

\begin{example} \label{ex:constcurv}
Any Riemannian manifold of constant curvature $C$ is IP. Its
curvature tensor $R^C$ has rank $2$ when $C \ne 0$.
\end{example}

\begin{example}\cite{GLS} \label{ex:phi}
Let $\phi$ be a linear isometry of $\Rn$ with $\phi^2 =\id$
(all the eigenvalues of such a $\phi$ must be $\pm 1$),
and let $C \ne 0$. Then an algebraic curvature tensor $R^C_{\phi}$ defined by
$R^C_{\phi}(X, Y) = R^C (\phi X, \phi Y)$ is IP, and $\rk \,R^C_{\phi} = 2$.
\end{example}

\begin{example}\cite{GLS, IP} \label{ex:f(t)}
A Riemannian manifold $M^n$ with a metric of a warped product
\begin{equation}\label{eq:f(t)}
ds^2 = dt^2 + f(t) ds_K^2,
\end{equation}
where $ds_K^2$ is a metric of constant curvature $K$
and $f(t) = Kt^2 + At + B > 0$, is IP. Its curvature tensor has the form
$R^{C(t)}_{\phi}$, with $C(t) =(4KB - A^2)/(4 f(t)^2)$.
For every point $x \in M^n, \; \phi$ is a reflection of the tangent space $T_xM^n$ in
the hyperplane orthogonal to $\partial/\partial t$.
\end{example}

In Example~\ref{ex:f(t)}, all but one eigenvalues of $\phi$ are $+1$. Clearly,
if all the eigenvalues of $\phi$ are the same ($\phi = \pm \id$), the resulting
algebraic curvature tensor (or manifold) has constant curvature. On the other hand,
no IP curvature tensors $R^C_{\phi}$ of Example~\ref{ex:phi}, with $\phi$ having
more than one eigenvalue $+1$ and more than one eigenvalue $-1$, can locally be
the curvature tensor of a Riemannian manifold \cite{GLS}.

Note that the metric \eqref{eq:f(t)} is not of constant curvature, unless
$4 KB - A^2 = 0$, but is conformally flat.

The IP manifolds were introduced and classified in dimension $4$
by Ivanov, Petrova \cite{IP} (hence the name).
Shortly after, in \cite{GLS}, Gilkey, Leahy and
Sadofsky using powerful topological methods classified all the IP algebraic curvature
tensors and manifolds of dimensions $n \ge 9$  and $n =5, 6$. Later, in \cite{G2},
Gilkey extended the result of \cite{GLS} to $n=8$, and gave a detailed description
of all possible eigenvalue structures of $R(\pi)$ when $n = 7$. The case $n =7$ was
further studied in \cite{GS} using spinors.

In this paper, we complete the case $n =7$:

\begin{theorem*} Any nonzero IP algebraic curvature tensor in $\RS$ has rank $2$.
\end{theorem*}

This, combined with the results of \cite{GLS, G2, IP}, gives the following
classification:

\begin{corollary*} {\ }
\begin{enumerate}
    \item Any nonzero IP algebraic curvature tensor in $\Rn,\; n \ne 4$, has rank $2$ and is of the
form $R^C_{\phi}$ of Example~\ref{ex:phi}.
    \item Any Riemannian IP manifold $M^n, \; n \ge 4$, is either of constant curvature, or is
locally isometric to the warped product \eqref{eq:f(t)}.
\end{enumerate}
\end{corollary*}

Note that the case $n=2$ is of no interest: any algebraic curvature tensor (any
Riemannian manifold) is IP. In dimension $3$, IP algebraic curvature tensors can be easily
classified (see \cite[Remark 1]{IP}): they are either of constant curvature, or those whose
Ricci tensor has rank $1$ (this fits the construction of Example~\ref{ex:phi}).
However, the class of IP Riemannian manifolds of dimension $3$ with such a Ricci tensor
is much wider than in Example~\ref{ex:f(t)} (see Section~\ref{s:dim3} for
discussion and some examples). In dimension $4$, there exist IP algebraic curvature tensors
of rank four (see \cite{IP} for classification, and \cite{GS} for construction using
spinors), but only those of rank $2$ can be realized as the curvature tensors of
$4$-dimensional manifolds
\footnote{The author is thankful to P.Gilkey for pointing out that the Corollary
was not stated correctly in the original version.}.

The IP algebraic curvature tensors were also extensively studied in pseudoriemannian and in
complex settings. We refer to \cite{Gb, GI1, GI2, GZ} for results in these directions.

The proof of the Theorem is given in Section~\ref{s:proof}. In Section~\ref{s:dim3}, we
study three-dimensional IP manifolds.

\section{Proof of the Theorem}\label{s:proof}

Let $R$ be an IP algebraic curvature tensor in $\RS$, whose rank is bigger than $2$. For
every two-plane $\pi \in G^+(2,7)$, the symmetric operator $R(\pi)^2$ has an odd-dimensional
kernel and some negative eigenvalues, $-\la_j^2$, each of an even multiplicity
$n_j, j = 1, \ldots, p$. Let $E_j(\pi)$ be the corresponding eigenspaces, with $\dim E_j(\pi) = n_j$.
Label the $n_j$'s in a non-decreasing order and call the ordered set
$(n_0 = \dim\Ker\,R(\pi), n_1,\ldots, n_p)$ the \emph{eigenvalue structure} for $R$.
Then, according to \cite[Theorem 0.4, 1a), (3)]{G2}, one has only two possibilities:
\begin{enumerate}
\item[(a)] the kernel is one-dimensional and $n_1 = 2, \, n_2 = 4$: the eigenvalue structure $(1,2,4)$;
\item[(b)] the kernel is three-dimensional and $n_1 = 4$: the eigenvalue structure $(3,4)$.
\end{enumerate}
We want to show that no IP algebraic curvature tensor with such eigenvalue structures can exist.
We will assume, in the both cases, that the eigenvalue of $R(\pi)^2$ of multiplicity $4$ is $-1$.
In an appropriate orthonormal basis for $\RS$, the matrix of the operator $R(\pi)$ has the following
normal form, respectively:
\begin{equation}\label{eq:normalform}
\mathrm{(a)}
\left(\begin{matrix} \mathcal{J}&0& 0\\0& \alpha J &0 \\0&0&0 \end{matrix}\right), \qquad
\mathrm{(b)}
\left(\begin{matrix} \mathcal{J}&0\\0&0\end{matrix}\right),
\quad \text{where} \;
\mathcal{J} = \left(\begin{matrix} 0&I_2\\-I_2&0\end{matrix}\right),\;
J = \left(\begin{matrix} 0&1\\-1&0\end{matrix}\right), \; \alpha \ne 0, \pm 1.
\end{equation}
For arbitrary $X, Y$, the normal form of the matrix of $R(X, Y)$ is the one above multiplied by
$\|X \wedge Y\|$.

The proof goes as follows. We start with the case (a) (Section~\ref{ss:124}). The first step
(Lemma~\ref{l:linker}) is to show that the kernel of $R(X,Y)$ is spanned by a vector depending
linearly on $X$ and $Y$. Next, in Lemma~\ref{l:existsU}, we prove that for any nonzero vector $X$,
the set $\K(X) = \cup_{\pi \ni X}\Ker \, R(\pi)$ is a linear space of dimension six. What is more, there
exists an orthogonal operator $U$ in $\RS$ such that $UX$ is a normal vector to $\K(X)$, for all
$X \ne 0$.
The key step in the proof is Lemma~\ref{l:const} saying that, for any
nonzero $X$ and any two-plane $\pi \ni X$, the two-dimensional eigenspace $E_1(\pi)$ of $R(\pi)$
contains the vector $UX$. It follows that $E_1(\pi) = U(\pi)$. We then show that
$U$ is symmetric, and that the tensor $R$ splits on two: $R^{\pm\alpha}_U$ (as in Example~\ref{ex:phi}),
and the remaining part, which is an IP algebraic curvature tensor with the eigenvalue structure
$(3, 4)$, hence reducing (a) to (b).

The case (b) is done in Section~\ref{ss:34} by a brute force of linear algebra, using the fact
that for all $X, Y$, the operator $R(X,Y)$ satisfies $R(X,Y)^3 + \W R(X,Y)= 0$, which follows
from \eqref{eq:normalform}.

\subsection{Case (a), the eigenvalue structure $(1,2,4)$} \label{ss:124}

We start with a brief introduction from commutative algebra. Let $\mathbf{D}$ be an integral domain (an
associative commutative ring with a $1$ and without zero divisors). A noninvertible element
$p \in \mathbf{D}$ is prime, if it  generates a prime ideal
($p \mid ab \implies p \mid a$ or $p \mid b$), and is irreducible, if $p = ab$ implies that
either $a$ or $b$ is invertible. The domain $\mathbf{D}$ is a \emph{unique factorization domain}
(UFD), if all irreducibles are primes and every element of $\mathbf{D}$ is a finite product
of irreducibles. In a UFD, every element $a$ can be represented in the form
$a = u \prod_i p_i^{m_i}$, with $p_i$ primes, $p_i \nmid p_j$, and $u$ invertible, and
such a representation is unique up to invertible elements. In particular, in a UFD, there
defined (up to invertibles) the greatest common divisor of a finite set of elements.
Also, for any four elements
$a_{11}, a_{12}, a_{21}, a_{22}$ satisfying $a_{11}a_{22} = a_{12}a_{21}$, there exist
$b_1, b_2, c_1, c_2$ such that $a_{ij} = b_ic_j$. Inductively, this implies the
following fact ($\mathbf{D}^n$ is a free module of rank $n$ over $\mathbf{D}$):

\begin{fact} \label{f:matrix}
Let $W$ be $n \times n$ matrix of rank $1$ \emph{(}all the $2 \times 2$ minors vanish\emph{)}
over a \emph{UFD} $\mathbf{D}$. Then
there exist $a, b \in \mathbf{D}^n$ such that $W = ab^t$.
If, in addition, $W$ is symmetric, then there exist $a \in \mathbf{D}^n, \; r \in \mathbf{D}$
such that $W = raa^t$.
\end{fact}

We will use the fact that a polynomial ring over reals is a UFD and the
Nagata Theorem\cite{Nag}:

\begin{fact}\label{f:nagata}
The ring $\mathbf{R} = \mathbb{R}[x_1, \ldots, x_n]/(\sum_{i=1}^n x_i^2)$
is a \emph{UFD}, when $n \ge 5$.
\end{fact}

Back to IP algebraic curvature tensors, we start by proving that the kernel of
$R(X,Y)$ depends linearly on $X$ and $Y$. More precisely:

\begin{lemma}\label{l:linker}
There exists a bilinear skew-symmetric map
$B:\RS \times \RS \to \RS$ such that
\begin{align}
\|B(X,Y)\|^2 &= \|X \wedge Y\|^2, &\quad &\text{for all $X, Y \in \RS$}
\label{eq:bxy^2} \\
\Ker \,R(X,Y) &= \Span (B(X,Y)), &\quad &\text{for all $X \nparallel Y \in \RS$}.
\label{eq:kerR}
\end{align}
\end{lemma}

\begin{proof} For every pair of vectors
$X, Y \in \RS$, define a symmetric operator $W(X,Y): \RS \to \RS$ by
$$
W(X, Y)= (R(X,Y)^2 + \|X \wedge Y\|^2)(R(X,Y)^2 + \a^2\|X \wedge Y\|^2).
$$
For arbitrary nonparallel $X, Y$, the operator $W(X, Y)$ has rank $1$,
with a nonzero eigenvalue $\a^2 \|X \wedge Y\|^4$ (this follows from
\eqref{eq:normalform}). The corresponding eigenvector spans $\Ker \, R(X, Y)$.

The matrix of $W(X, Y)$ can be viewed as a matrix over the ring
$\mathbf{K}=\mathbb R[x_1, \ldots, x_7, y_1, \ldots, y_7,]$ of polynomials in
$14$ variables, the coordinates of $X$ and $Y$
(all its entries are homogeneous polynomials of degree $4$).
By Fact~\ref{f:matrix},
there exist a polynomial $f(X,Y)$ and a $7$-vector $P(X, Y)$ with
polynomial components such that $W(X, Y) = f(X,Y) P(X,Y) P(X,Y)^t$.
For any nonparallel $X$ and $Y$, the vector $P(X,Y)$ spans
$\Ker \, R(X, Y)$. As $\Tr\, W(X, Y) = \a^2\|X \wedge Y\|^4$ and the
polynomial $\|X \wedge Y\|^2$ is irreducible in $\mathbf{K}$,
we have two possibilities for $f$ (up to multiplication by a positive constant):
either $f = 1$, or $f =\|X \wedge Y\|^2$ (the case $f = \|X \wedge Y\|^4$
is not possible, as then the vector $P(X, Y)$, which spans $\Ker \, R(X, Y)$,
is constant. But if $Z \in \RS$ is in the kernel of all the $R(X,Y)$'s,
then $R(Z,\cdot)$ is zero).

We want to show that the case $f =1$ leads to a contradiction. Assume $f=1$,
hence $W(X, Y) = P(X,Y) P(X,Y)^t$ for a polynomial vector $P \in \mathbf{K}^7$.
As the entries of $W$ are homogeneous in $X$, of degree $4$, and homogeneous
in $Y$, of degree $4$, the $P_i$'s, the components of $P$, must be polynomials
homogeneous in $X$, of degree $2$, and homogeneous in $Y$, of degree $2$
(each component of $P$ is a linear combination of terms $x_ix_jy_ky_l$).
We also have $\sum_{i=1}^7 P_i^2(X, Y) = \Tr \,W(X, Y) = \a^2 \|X \wedge Y\|^4$.

For every nonzero $X \in \RS$, define the the subset $\K(X) \subset \RS$ as follows:
$$
\K(X) = \cup_{Y \nparallel X} \Ker\,  R(X,Y) = \cup_{\pi \ni X} \Ker\,  R(\pi) =
\{\Span(P(X, Y)) : Y \perp X,\, \|Y\|=1\}.
$$
The set $\K(X)$ is a cone over the image of the sphere $S^5 \subset \mathbb R^6$
under the polynomial map, hence its complement $\RS \setminus \K(X)$ is open and dense.
It follows that the set of pairs $(X,Z) \in \RS \times \RS$ such that
$Z \notin \K(X),\; X \notin \K(Z)$ is nonempty (even dense). Let $(X, Z)$
be one such pair, and $S^5$ be a unit sphere in $X^\perp$.
Consider a map $V: S^5 \to X^\perp$ defined by $V(Y) = R(Z, P(X,Y))X$.
We have $\<Y, V(Y)\> = R(Z, P(X,Y),Y, X) = 0$, as $P(X, Y)$ is in the kernel
of $R(X,Y)$. Furthermore, $\Range (V) \not\ni 0$. Indeed, if $V(Y_0) = 0$
for some $Y_0 \in S^5$, then for all $T \in \RS, \; 0 = R(Z, P(X,Y_0), X, T)$.
As $P(X,Y_0) \ne 0$ (since $\|P(X,Y_0)\|^2 = \a^2\|X \wedge Y_0\|^4 = \a^2\|X\|^4$),
and $P(X,Y_0) \nparallel Z$ (since $Z \notin \K(X)$), this implies
that $X \in \K(Z)$, which contradicts the choice of the pair $(X, Z)$. Now the map
$\hat V: S^5 \to S^5$ defined by $\hat V(Y)= V(Y)/\|V(Y)\|$ is even
($\hat V(-Y) = \hat V(Y)$), as such is $P(X, Y)$, and $\hat V(Y) \perp Y$.
This is not possible, since otherwise the homotopy $Y \cos t + \hat V(Y)\sin t$ joins
the identity map of $S^5$ with the one of an even degree.

It follows that $f = \|X \wedge Y\|^2$, and so
$W(X, Y) = \|X \wedge Y\|^2 P(X,Y) P(X,Y)^t$. Comparing the degrees, we find
that all the components of the polynomial vector $P(X, Y)$ are linear in $X$
and in $Y$, so each $P_i(X, Y)$ is a bilinear form on $\RS$. Also,
$\sum_{i=1}^7 P_i^2(X, Y) = \|X \wedge Y\|^{-2} \Tr \,W(X, Y) = \a^2\|X \wedge Y\|^2$,
which implies $P(X,X)=0$, so $P$ is skew-symmetric. Finally, for $X \nparallel Y$,
$P(X, Y)$ is a nonzero eigenvector of $W(X, Y)$, hence it spans the kernel of
$R(X,Y)$.

Now define the map $B$ by setting $B(X, Y)=\a^{-1} P(X, Y)$.
\end{proof}

\begin{lemma}\label{l:existsU}
1) For $X \ne 0$, the set $\K(X) = \cup_{\pi \ni X} \Ker \, R(\pi)$ is
a six-dimensional subspace of $\RS$.\newline
2) There exists an orthogonal operator $U$ on $\RS$ such that for all
$X \ne 0$, the vector $UX$ is orthogonal to $\K(X)$, or equivalently,
\begin{equation}
UX \perp \Ker \, R(X, Y) \quad \text{for all}\; Y \nparallel X.
\end{equation}
\end{lemma}
\begin{proof} Let $\Ro$ be an orthogonal sum of $\mathbb R e_0$ and $\RS$,
with $p: \Ro \to \RS$ the orthogonal projection.
Define a bilinear map
$\overline{B}: \RS\times \RS \to \Ro$ as follows:
$$
\overline{B}(X, Y)= \<X,Y\>e_0 + B(X, Y),
$$
where $B$ is the map from Lemma~\ref{l:linker}. Then
for all $X, Y \in \RS$ and all $Y_1 \perp Y_2$,
\begin{equation}\label{eq:barB}
\|\overline{B}(X,Y)\|^2 = \|X\|^2 \| Y\|^2, \quad
\<\overline{B}(X,Y_1), \overline{B}(X,Y_2)\> = 0,
\end{equation}
(the first equation follows from \eqref{eq:bxy^2}, the second one follows
from the first one), so $\overline{B}$ is a normed bilinear map.
For every $X \in \RS$, define an operator $A_X: \Ro \to \RS$ by
$$
\<A_XZ, Y\> = \<\overline{B}(X,Y), Z\>,
$$
where $Z \in \Ro, \, Y \in \RS$. Then from \eqref{eq:barB}, for all
$X$ and all $X_1 \perp X_2$,
\begin{equation}\label{eq:AA^t}
A_XA_X^t = \|X\|^2 \,\id_{\RS}, \quad A_{X_1}A_{X_2}^t \in \mathfrak o(7),
\end{equation}
where $\mathfrak o(7)$ is the linear space of skew-symmetric operators in $\RS$.
In particular, from the first equation, $\rk \, A_X = 7$ when $X
\ne 0$. This proves assertion 1 of the Lemma. Indeed, by
Lemma~\ref{l:linker}, $\K(X) = \cup_{Y \nparallel X} \Ker \, R(X, Y) =
\Range\, B(X, \cdot)$, which is a linear subspace of $\RS$, of
dimension at most $6$ (as $B$ is skew-symmetric). If a nonzero
vector $Z' \in \RS$ is orthogonal to this subspace, then the
vector $Z = 0e_0 + Z' \in \Ro$ is in the kernel of $A_X$, which is
of dimension $1$. Thus $\dim \K(X) = 6$.

Fix a unit vector $X_0 \in \RS$. The kernel of the operator $A_0 = A_{X_0}$
is one-dimensional. Let $Z_0$ be a unit vector in this kernel. For a nonzero
vector $X \perp X_0$, let $Y \in \RS$ be a nonzero vector from the kernel of the
skew-symmetric operator $A_0A_X^t$, so that $A_0A_X^tY = 0$. As $\rk\, A_0 = 7$,
with $\Ker \,A_0$ spanned by $Z_0$, it follows that $A_X^tY$ is parallel to $Z_0$.
Since $Y \ne 0$ and $\rk \, A_X^t = 7$, the vector $A_X^tY$ is nonzero, and up to
scaling, we can choose $Y$ in such a way that $A_X^t Y = \|X\|^2 Z_0$. Acting on the
both sides by $A_X$ we find $Y = A_XZ_0$, so for all nonzero $X \perp X_0$,
\begin{equation}\label{eq:A_Xetc}
A_XA_0^tA_XZ_0 = -A_0A_X^tA_XZ_0 = 0, \quad A_XZ_0 \ne 0.
\end{equation}
Define a linear operator $V: \RS \to \Ro$ by
\begin{equation}\label{eq:defU}
VX_0 := Z_0, \quad VX := -A_0^t A_X Z_0 \quad \text{for}\; X \perp X_0.
\end{equation}
We want to show that $U = p \circ V : \RS \to \RS$ is the sought orthogonal operator.
The operator $V$ has the following properties:

(i) For all $X \in \RS, \quad VX \in \Ker\, A_X$.

(ii) $V$ is orthogonal (note that it acts between Euclidean spaces of different
dimension):
$$
\|VX_0\| = 1, \quad \text{and for} \; X \perp X_0: \quad VX \perp VX_0, \; \|VX\| = \|X\|.
$$

To check (i), take an arbitrary $X \perp X_0$ and $t \in \mathbb R$. We have:
$$
A_{tX_0 +X}V(tX_0 +X) = (tA_0 + A_X)(tZ_0 - A_0^t A_X Z_0)=
t^2 A_0Z_0 +t(A_XZ_0-A_0A_0^t A_X Z_0) - A_XA_0^t A_X Z_0= 0,
$$
by \eqref{eq:A_Xetc} and \eqref{eq:AA^t}.

The first equation of (ii) immediately follows from \eqref{eq:defU}. For the second one, we have:
$\<VX, VX_0\> = \<-A_0^t A_X Z_0, Z_0\> = - \<A_X Z_0, A_0 Z_0\> = 0$.
To check the third one, consider the vector $A_X^tA_XZ_0$. As
$A_0(A_X^tA_XZ_0) = 0$ by \eqref{eq:A_Xetc}, $A_X^tA_XZ_0 = f(X) Z_0$
for some function $f$. Acting on both sides by $A_X$ we get
$(f(X) - \|X\|^2)A_XZ_0 = 0$, so $A_X^tA_XZ_0 = \|X\|^2 Z_0$, since
$A_XZ_0 \ne 0$ by \eqref{eq:A_Xetc}. Then
$$
\|VX\|^2 = \<A_0^t A_X Z_0, A_0^t A_X Z_0\> = \<A_X Z_0, A_0 A_0^t A_X Z_0\> =
\<A_X Z_0, A_X Z_0\> = \<A_X^tA_XZ_0, Z_0\> = \|X\|^2,
$$
as required.

From property (i) it follows that for all
$Y \in \RS,\, 0 = \<A_XVX,Y\>= \<\overline{B}(X,Y), VX\>$. Taking $Y=X$ we get
$VX \perp e_0$, for all $X$ (as $B$ is skew-symmetric). So $\Range \, V = \RS$,
and the operator $U:\RS \to \RS$ defined by $U = p \circ V$ is orthogonal
($U$ acts exactly as $V$, but with a different codomain).
Moreover, as $0 = \<\overline{B}(X,Y), VX\> = \<B(X,Y), UX\>$, we have
$UX \perp B(X, Y) = \Ker \, R(X, Y)$, for all
$X, Y$.
\end{proof}

\begin{remark} From the proof of Lemma~\ref{l:existsU}, it is easy to see that the map
$\Phi: \RS \to \Hom(\Ro, \Ro)$ defined by $\Phi(X)Z = A_XZ + \<VX,Z\>e_0$
has the property $\Phi(X)\Phi(X)^t = \|X\|^2 \id_{\Ro}$, and so
the map $\phi(X) = \Phi(X_0)^t\Phi(X)$ defined on the six-space $X_0^\perp$
satisfies $\phi(X)^2 = -\|X\|^2 \id_{\Ro}$. Thus $\phi$ can be extended to
a 
representation of the Clifford algebra $\Cl(6)$ in $\Ro$, which is a
restriction of that for the Clifford algebra $\Cl(7)$, which, in turn, is
equivalent to the right (or to the left) multiplication by imaginary
octonions in the octonion algebra $\mathbb{O}$. One can then show that,
identifying $\RS$ with the space of imaginary octonions,
$B(X, Y)$ is the imaginary part of $XY$, up to orthogonal transformations.
\end{remark}

For every pair of nonparallel vectors $X, Y$, let $E_1(X,Y)$ be the
two-dimensional eigenspace of $R(X, Y)^2$ with the eigenvalue $-\a^2 \W$
(that is, $E_1(X,Y) = E_1(\pi)$, where $\pi = \Span(X, Y)$).

\begin{lemma}\label{l:const} $E_1(X,Y)= \Span (UX, UY)$,
where $U$ is the orthogonal operator introduced in Lemma~\ref{l:existsU}.
\end{lemma}
\begin{proof}
Fix a unit vector $X$. Introduce a new variable $t$, and define,
for every $Y \in X^\perp$, the operators $G(Y), M(Y, t):\RS \to \RS$ by
\begin{gather*}
G(Y) = R(X,Y)^2+\|Y\|^2 \id -B(X,Y)B(X,Y)^t, \\
M(Y,t)=(R(X,Y) + t \, \a \, \id )\,G(Y)= R(X,Y)\,G(Y) + t\, \a \, G(Y),
\end{gather*}
where $B$ is the map from Lemma~\ref{l:linker} spanning the kernel of $R(X,Y)$.
Note that the operator $G(Y)$ is symmetric, while the operator
$R(X,Y)\,G(Y)$ is skew-symmetric.

At this point, it will be more convenient to switch from operators to matrices
fixing some orthonormal basis for $\RS$. With a slight abuse of language, we will
use the same notation for an operator and its matrix. For $Y \in X^\perp$, let
$y_1, \ldots, y_6$ be its coordinates with respect to an orthonormal basis for
$X^\perp$ (which is not related to the chosen orthonormal basis for $\RS$).
Denote
$\mathbb R[Y]=\mathbb R[y_1, \ldots, y_6]$ and $\mathbb R[Y,t]=\mathbb R[y_1, \ldots, y_6,t]$
the corresponding polynomial rings.

From definition, it is clear that all the entries of $G(Y)$ and $R(X,Y)\,G(Y)$ are
homogeneous polynomials of the $y_i$'s, of degree $2$ and $3$, respectively.

From \eqref{eq:normalform}, the normal forms of the matrices
$G(Y), R(X,Y)\,G(Y)$ and $M(Y,t)$ are, respectively,
\begin{equation}\label{eq:similar}
\|Y\|^2 (1-\a^2) \left(\begin{matrix} 0&0& 0\\0& I_2 &0 \\0&0&0 \end{matrix}\right), \qquad
f(Y)\|Y\| \left(\begin{matrix} 0&0& 0\\0& J &0 \\0&0&0 \end{matrix}\right), \qquad
f(Y) \left(\begin{matrix} 0&0& 0\\0& t I_2 + \|Y\| J &0 \\0&0&0 \end{matrix}\right),
\end{equation}
in the same basis, where
$f(Y) = \|Y\|^2 \a(1-\a^2),\; J = \left(\begin{smallmatrix} 0&1\\-1&0\end{smallmatrix}\right)$,
and the nonzero blocks are in the $5$th and the $6$th rows and columns.

As it follows from \eqref{eq:similar}, for every nonzero $Y \perp X$,
both $G(Y)$ and $R(X,Y)\,G(Y)$ have rank $2$. Moreover, the two-space
$E_1(X,Y)$ is the image of the operator $R(X,Y)\,G(Y)$ and is the
eigenspace of the operator $G(Y)$, with the eigenvalue $(1-\a^2)\|Y\|^2$.

It follows from \eqref{eq:similar} that
\begin{equation} \label{eq:M^2}
M(Y,t)^2 -2 \a(1-\a^2)\|Y\|^2 t M(Y,t) = -\a^3(1-\a^2)\|Y\|^2(\|Y\|^2+t^2)G(Y).
\end{equation}
Moreover, as $M(Y,t)$ still has the normal form \eqref{eq:similar} for real $Y$ and
\emph{complex} $t$, the rank of the complex matrix $M(Y, i\|Y\|)$ is $1$ for all nonzero
$Y \in \RS$. So all the $2\times 2$ minors of the polynomial matrix $M(Y,t)$ vanish for
$t = i\|Y\|$. Any such minor has a form $q(Y,t)=f(Y) + t g(Y) + t^2 h(Y)$, with $f, g, h$
real polynomials. As $q(Y, i\|Y\|) = 0$, we get $g(Y) = 0,\, f(Y) = \|Y\|^2 h(Y)$, hence
every $2\times 2$ minor of $M(Y,t)$ is divisible by $t^2 + \|Y\|^2$ in the polynomial
ring $\mathbb R[Y,t]$.

Let $\mathbf I \subset \mathbb R[Y,t]$ be the ideal generated by $t^2 + \|Y\|^2$, and
let $\RR = \mathbb R[Y,t] / \mathbf I$, with $\pi: \mathbb R[Y,t] \to \RR$ the
natural projection. Note that for every element $a \in \RR$, there is a unique
pair of polynomials $p, q \in \mathbb R[Y]$ such that $\pi (p + t q) = a$.

Consider the $7 \times 7$ matrix $\M = \pi(M)$, with entries from $\RR$. As all
the $2\times 2$ minors of $M(Y,t)$ are in $\mathbf I$, the rank of the matrix
$\M$ is $1$ ($\M$ is nonzero, since nonzero entries of $M$ are at most linear
in $t$). Projecting the equation \eqref{eq:M^2} to $\RR$, we obtain
\begin{equation}\label{eq:m^2}
  \M^2= -2\a(1-\a^2)\,{\bar t\,}^3\,\M,
\end{equation}
where $\bar t = \pi(t)$.

By Fact~\ref{f:nagata}, the ring $\RR$ is a unique factorization domain. Let
$d \in \RR$ be the greatest common divisor of the entries of $\M$, and $\M = d \mathcal{L}$,
with the greatest common divisor of the entries of $\mathcal L$ being $1$. Let $L_1, L_2$ be
matrices with entries from $\mathbb R[Y]$ such that $\pi (L_1 +tL_2) = \mathcal L$.

From \eqref{eq:m^2}, $d \mid \bar t^3$, and so (as $\bar t$ is prime in $\RR$), $d = \bar t^m$,
where $m = 0, 1, 2, 3$. Consider these cases separately.

First show that $m > 0$. As $\rk \, \M = 1$, by Fact \ref{f:matrix}, there exist
$a,b \in \mathbf{R}^7$ such that
$\M = ab^t$. Reducing $M(Y,t) + M(Y,t)^t = 2 \a t G(Y)$ modulo $\mathbf I$ we get
$ab^t + ba^t = 2 \a \bar t \pi(G(Y))$. So for all $i, j = 1, \ldots, 7$,
\begin{equation}\label{eq:ab'+ba'}
  a_ib_j + b_ia_j = 2 \a \bar t \, \pi(G_{ij}(Y)).
\end{equation}
Taking $j=i$ in \eqref{eq:ab'+ba'} we find that $\bar t \mid a_ib_i$, so for every
$i$, at least one of $a_i, b_i$ is divisible by $\bar t$. If for some $i \ne j$,
$\bar t \nmid a_i, b_j$, then $\bar t \mid b_i, a_j$, and we come to a
contradiction with \eqref{eq:ab'+ba'}. It follows that either all the $a_i$'s,
or all the $b_i$'s are divisible by $\bar t$. In both cases, all the entries
of the matrix $\M = ab^t$ are divisible by $\bar t$, so $m >0$.

Assume that $m =3$. Lifting the equation $\M = \bar t ^3 \mathcal{L}$ to $\mathbb R[Y,t]$
we get, modulo $\mathbf I, \quad M(Y,t) = t^3 (L_1 + t L_2) = - t \|Y\|^2 (L_1 + t L_2)
= \|Y\|^2(\|Y\|^2L_2-tL_1)$, hence $M(Y,t) = \|Y\|^2(\|Y\|^2L_2-tL_1) + (t^2 + \|Y\|^2)\hat M$
for some matrix $\hat M$ with entries in $\mathbb R[Y,t]$. As
$M(Y,t)= R(X,Y) G(Y) + t \a G(Y)$, with all the entries of $R(X,Y) G(Y)$ and $G(Y)$
homogeneous polynomials of $Y$ of degree $3$ and $2$, respectively, we get
$R(X,Y) G(Y) = 0$, which contradicts the fact that $\rk \, R(X,Y) G(Y) = 2$, when
$Y$ is nonzero.

Let $m =2$. We have $M(Y,t) = t^2 (L_1 + t L_2) \mod \mathbf I$, so
$M(Y,t) = -\|Y\|^2 L_1  -\|Y\|^2 t L_2 + (t^2 + \|Y\|^2)\hat M = R(X,Y) G(Y) + t \a G(Y)$.
It follows that $G(Y) = \|Y\|^2 G_0$ for some constant symmetric matrix $G_0$ of rank
$2$. But then for all nonzero $Y \perp X$, the eigenspace $E_1(X,Y)$ is the same: it is
the eigenspace of the fixed matrix $G_0$. This contradicts
Lemma~\ref{l:existsU}: the set $\K(X) = \cup_{Y \perp X, Y \ne 0} \Ker \, R(X, Y)$ is
a six-dimensional subspace of $\RS$, hence for some $Y$, the subspaces $\Ker \, R(X, Y)$ and
$E_1(X,Y)$ have a nonzero intersection.

Finally, consider the case $m=1$. Then $M(Y,t) = t (L_1 + t L_2) \mod \mathbf I$, so
$M(Y,t) = -\|Y\|^2 L_2 + t L_1   + (t^2 + \|Y\|^2)\hat M$ for some matrices $L_1, L_2$,
with entries from $\mathbb R[Y]$, and a matrix $\hat M$, with entries from $\mathbb R[Y,t]$.
As $M(Y,t) = R(X,Y) G(Y) + t \a G(Y)$, it follows that $R(X,Y) G(Y) = -\|Y\|^2 L_2$,
and $L_2 = L_2(Y)$ must be a skew-symmetric matrix, of rank $2$ (when $Y \ne 0$),
whose entries are linear in $Y$.

We get a linear map $L_2$ from $\mathbb R^6 = X^\perp$ to $\mathfrak o_2(7)$, the set
of skew-symmetric $7 \times 7$ matrices of rank less than or equal to two. The map
$L_2$ is injective (as $\rk \, L_2(Y) = 2$ for all $Y \ne 0$), so by Lemma 2.2 of
\cite{GLS}, there exists a unit vector $\xi \in \RS$ such that
$L_2(Y)Z = \<\xi, Z\> L_2(Y)\xi -\<L_2(Y)\xi, Z\> \xi$, for all
$Z \in \RS$. In particular, taking $Z = L_2(Y)\xi$ we get
$L_2(Y)(L_2(Y)\xi) = -\|L_2(Y)\xi\|^2 \xi$. As $L_2(Y)\xi \ne 0$, unless $Y=0$
(otherwise $L_2(Y) = 0$), we find that, for all $Y \ne 0$,
$$
\xi \in \Range \,L_2(Y) = \Range \,R(X,Y) G(Y) = E_1(X, Y).
$$
For every nonparallel $X, Y,\quad \Ker\, R(X, Y) \perp E_1(X,Y)$ (they are
the eigenspaces of the symmetric operator $R(X, Y)^2$, with different
eigenvalues). It follows that $\xi \perp \cup_{Y\nparallel X} \Ker\, R(X, Y)= \K(X)$.
By Lemma~\ref{l:existsU}, $\K(X)$ is a six-dimensional subspace of $\RS$, whose
orthogonal complement is spanned by the vector $UX$.

It follows that $UX \in E_1(X, Y)$. Similarly, $UY \in E_1(X, Y)$. As the operator
$U$ is orthogonal (and, in particularly, nonsingular), the two-dimensional spaces
$\Span(UX, UY)$ and $E_1(X, Y)$ must coincide.
\end{proof}

Let $X, Y$ be any two orthonormal vectors. From Lemma~\ref{l:const} it follows that
$E_1(X,Y) = \Span(UX, UY)$. Then
$R(X, Y) UY = \varepsilon \a UX,\; \varepsilon = \pm 1$, and by continuity,
$\varepsilon$ is the same for all $X, Y$. Hence for any $Z$,
\begin{equation}\label{eq:Rxyuy}
R(X, Y, UY, Z) = -\varepsilon \a \<UX, Z\>.
\end{equation}
We claim that the operator $U$ is not only orthogonal, but also symmetric.
Indeed, an orthogonal operator in an odd-dimensional space has at least one
eigenvalue $\pm 1$. Replacing $U$ by $-U$, if necessary, we can assume that
there exists a unit vector $Y \in \RS$ such that $UY = Y$. Note that the
space $Y^\perp$ is an invariant subspace of $U$. The equation \eqref{eq:Rxyuy},
with $X, Z \in Y^\perp$ gives:
$$
R(X, Y, Y, Z) = -\varepsilon \a \<UX, Z\>.
$$
The left-hand side is symmetric with respect to $X, Z$, and so such is the
right-hand side. It follows that the operator $U$ is symmetric on its invariant
subspace $Y^\perp$, hence is symmetric on the whole $\RS$.

As $U$ is orthogonal and symmetric, $U^2 = \id$. Let now $R^\a_U$ be an
algebraic curvature tensor constructed as in Example~\ref{ex:phi}, with the
operator $U$ and the constant $\a$:
$$
R^\a_U (X, Y) Z = \a \,(\<UY, Z\> UX - \<UX,Z\> UY).
$$
Define an algebraic curvature tensor $\bar R = R - \varepsilon R^\a_U$. We
claim that $\bar R$ is IP, with the eigenvalue structure $(3,4)$. Indeed, for any
two orthonormal vectors $X, Y$, we have:

1. If $Z \in \Ker \, R(X, Y)$, then $Z \perp UX, UY$ (by assertion 2 of Lemma~\ref{l:existsU}),
and so
$\bar R(X, Y) Z = - \varepsilon R^\a_U (X, Y) Z
= - \varepsilon\a (\<UY, Z\> UX -  \<UX,Z\> UY)= 0$, which implies
$Z \in \Ker \, \bar R(X, Y)$.

2. Let $Z = UY$. By \eqref{eq:Rxyuy}, $R(X, Y) UY = \varepsilon \a UX$. So
$\bar R(X, Y) UY = \varepsilon \a UX - \varepsilon \a (\<UY, UY\> UX - \<UX,UY\> UY)
= 0$, as $U$ is orthogonal. The same is true for $Z = UX$.
So $UX, UY \in \Ker \, \bar R(X, Y)$.

3. If $Z \in E_2(X,Y) = \left(\Span(\Ker \, R(X, Y), UX, UY)\right)^\perp$, then
$\bar R(X, Y) Z = R(X, Y)Z \in E_2(X, Y)$, So $E_2(X, Y)$ is an invariant subspace
of $\bar R(X, Y)$, and the restriction of $\bar R(X, Y)$ to $E_2(X, Y)$ has the same
eigenvalues as those of $R(X, Y)$, namely $\pm i$, both with multiplicity $2$.

\subsection{Case (b), the eigenvalue structure $(3,4)$} \label{ss:34}

Following \cite{G2}, for a nonzero $X \in \RS$, define a subset $\A(X) \in \RS$ by
$$
\A(X) = \cap_{Y \nparallel X} \Ker\,  R(X,Y).
$$
\begin{lemma}\label{l:u0}
There exists an an open, dense set $\mathcal{S} \subset \RS$ such that $\A(X) = 0$, when $X \in \mathcal{S}$.
\end{lemma}
\begin{proof}
Let for a given $X \ne 0, \; Z \in \A(X)$. Then for any $U, V \in \RS, \; R(X, U, V, Z) = 0$ and $R(X, V, Z, U) =
0$. So, by the first Bianchi identity, $R(X, Z, U, V) = 0$, that is, the operator $R(X, Z)$ is zero. It follows that
$Z \parallel X$. If $Z \ne 0$, then $R(X, V, X, U)=0$ for all $U, V \in \RS$, and, in particular, the curvature
on any two-plane in $\RS$ containing $X$ must vanish. The set of $X$'s with this property is closed. If
it has a nonempty interior, then the sectional curvature vanishes identically, and so $R=0$, which is a
contradiction.
\end{proof}
As it follows from \eqref{eq:normalform},  for all $X, Y \in \RS$,
\begin{equation}\label{eq:hj}
(R(X, Y))^3 + \W R(X,Y) = 0.
\end{equation}
Fix two orthonormal vectors $X, Y \in \RS$, with $X \in \mathcal{S}$, and choose
an orthonormal basis for $\RS$ in such a way that the matrix of $R(X,Y)$ is
$$
K = \left(\begin{matrix} \mathcal{J}&0\\0&0\end{matrix}\right),
\quad \text{where} \;
\mathcal{J} = \left(\begin{matrix} 0&I_2\\-I_2&0\end{matrix}\right).
$$
Let $\mathcal{W} = (\Span(X, Y))^\perp$. For a vector $Z \in \mathcal{W}$,
let $L(Z) = \left(\begin{smallmatrix}A(Z)&B(Z)\\-B(Z)^t&C(Z)\end{smallmatrix}\right)$
be the matrix of $R(X,Z)$, with $A(Z), C(Z)$ skew-symmetric
$4 \times 4$- and $3 \times 3$-matrices, respectively, and $B(Z)$ a
$4\times 3$-matrix, all depending linearly on $Z \in \mathcal{W}$.
The equation \eqref{eq:hj} gives
\begin{equation}\label{eq:hj2}
(yK +L(Z))^3 = -(y^2 + \|Z\|^2) (yK +L(Z)), \quad \text{for all $y \in \mathbb{R}$},
\end{equation}
Expanding \eqref{eq:hj2} by powers of $y$ we find:
\begin{gather}
K^2L(Z) + L(Z)K^2 + K L(Z) K = -L(Z), \label{eq:kl1}\\
L(Z)^2K + KL(Z)^2 + L(Z) K L(Z) = -\|Z\|^2K, \label{eq:kl2}\\
L(Z)^3 = - \|Z\|^2L(Z).\label{eq:kl3}
\end{gather}
From \eqref{eq:kl1} it follows that $C(Z)=0$ and $\mathcal{J}A(Z)\mathcal{J}=A(Z)$,
hence $A(Z) = \left(\begin{smallmatrix}a(Z)J&b(Z)J\\b(Z)J&-a(Z)J\end{smallmatrix}\right)$,
where $J = \left(\begin{smallmatrix} 0&1\\-1&0\end{smallmatrix}\right)$, and $a, b$ are
linear functionals on $\mathcal W$.
Then from \eqref{eq:kl2}, \eqref{eq:kl3},
\begin{gather}
B(Z)^t\mathcal{J}B(Z) = 0, \label{eq:bjb}\\
B(Z)B(Z)^t -\mathcal{J}B(Z) B(Z)^t\mathcal{J} = (\|Z\|^2-a(Z)^2-b(Z)^2)I_4 \label{eq:bb'}\\
B(Z)^tA(Z)B(Z)=0 \label{eq:bab}.
\end{gather}
Equation \eqref{eq:bjb} implies that the column space of the matrix
$B(Z)$ is an isotropic subspace of $\mathcal J$, hence $\rk \,B(Z) \le 2$.
If the set of vectors $Z \in \mathcal W$ with $\rk \,B(Z) < 2$ has a nonempty
interior, then $\rk\, B(Z) <2$ for all $Z$, and the matrix on the right-hand side
of \eqref{eq:bb'} has rank at most two. It follows that $a(Z)^2 + b(Z)^2 =\|Z\|^2$,
which is not possible for two linear functionals on a $5$-space. So, for an open,
dense set of vectors $Z \in\mathcal W,\; \rk\, B(Z) = 2$.

Multiplying the equation \eqref{eq:bb'} by $B(Z)$ from the right
and using \eqref{eq:bjb}, we get
$$
B(Z)B(Z)^tB(Z) = (\|Z\|^2-a(Z)^2-b(Z)^2)B(Z).
$$
This equation, together with the fact that $\rk\, B(Z) = 2$ for almost all
$Z \in\mathcal W$, implies that the singular numbers of the $4 \times 3$
matrix $B(Z)$ are $c(Z), c(Z), 0$, where $c(Z) =\sqrt{\|Z\|^2-a(Z)^2-b(Z)^2}$.

We need the following Lemma:

\begin{lemma}\label{l:cc0}
Let $\mathcal{V}$ be a linear space of $4 \times 3$-matrices $B$
whose singular numbers are $c, c, 0$ \emph{(}where $c =c(B) \ge 0$\emph{)}, and
such that
$\cap_{B \in \mathcal V} \Ker\, B = 0$. Then $\dim \mathcal{V} \le 3$.
\end{lemma}
In fact, up to orthogonal transformations, $\mathcal{V}$ is a subspace of
the space of $3 \times 3$ skew-symmetric matrices, with a zero row added
at the bottom.
\begin{proof}[Proof of Lemma~\ref{l:cc0}]
Let $B_1, B_2, B_3, B_4$ be linearly independent matrices in $\mathcal{V}$.
Up to orthogonal transformation and normalization, we can assume that
$$
B_1 =\left(\begin{matrix}I_2&0\\0&0\end{matrix}\right),
$$
and $\Tr \,B_1^tB_i =0$ for $i = 2, 3, 4$. The fact that $\rk (B_i + t B_1) \le 2$
implies that $B_i = \left(\begin{smallmatrix} Q_i&u_i\\T_i&0\end{smallmatrix}\right)$
for some $2 \times 2$-matrices $Q_i, T_i$ and $2$-vectors $u_i$ satisfying $T_iu_i =0$.
At least one of the $u_i$'s must be nonzero by assumption, and we can assume, up to
orthogonal transformation and up to taking appropriate linear combinations, that
$u_2 = (p, 0),\; p \ne 0$ and $u_4 = 0$. The fact that the singular numbers of the
matrix $(B_2+sB_4) + t B_1$ are $c, c, 0$ (for some $c$ depending on $t$ and $s$), together
with the condition $\Tr \,B_1^tB_i =0$, gives
$$
B_2 +s B_4 =\left(\begin{matrix}0&0&p\\0&0&0\\0&q_1(s)&0\\0&q_2(s)&0\end{matrix}\right),
$$
for some linear functions $q_1(s), q_2(s)$ satisfying $q_1(s)^2 + q_2(s)^2 = p^2$.
It follows that $q_1$ and $q_2$ are constants, hence $B_4 = 0$.
\end{proof}

Note that Lemma~\ref{l:cc0} applies in our situation, as
$\cap_{Z \in \mathcal W} \Ker\, B(Z) = 0$. Otherwise, if $u$ is a nonzero
vector with $B(Z)u = 0$ for all $Z \in \mathcal W$, then the set
$\A(X)$ contains a nonzero vector $(0,0,0,0,u)$, which contradicts the
choice of $X \in \mathcal{S}$.

By Lemma~\ref{l:cc0}, we can find two orthonormal vectors
$Z_1, Z_2 \in \mathcal{W}$ such that $B(Z_1) = B(Z_2) = 0$. It then follows
from \eqref{eq:bb'} that $a(Z)^2+b(Z)^2=\|Z\|^2$, for all
$Z \in \Span(Z_1, Z_2)$, so we can choose $Z_1, Z_2$ in such a way that
$a(Z_1)=b(Z_2)=1,\; a(Z_2)=b(Z_1)=0$.

Now for any $Z' \in \mathcal{W}$, the equation \eqref{eq:bab} with
$Z = Z' + t_1Z_1+t_2Z_2$ gives
$$
B(Z')^tA(Z_1)B(Z')=B(Z')^tA(Z_2)B(Z')=0.
$$
As a common isotropic subspace of the matrices
$$
A(Z_1)= \left(\begin{matrix}J&0\\0&-J\end{matrix}\right), \quad
A(Z_2)= \left(\begin{matrix}0&J\\J&0\end{matrix}\right), \quad
\mathcal{J} = \left(\begin{matrix} 0&I_2\\-I_2&0\end{matrix}\right)
$$
is at most one-dimensional, the latter equation, together with \eqref{eq:bjb},
implies that $\rk \,B(Z') \le 1$ for all $Z' \in \mathcal{W}$. This is a
contradiction with the fact that $\rk \,B(Z) =2$ for a generic $Z \in \mathcal{W}$.

\section{IP manifolds of dimension three} \label{s:dim3}

In the study of Riemannian IP manifolds, the three-dimensional case is
exceptional. In dimension $n \ge 4$, the most difficult part is algebraic.
Once all the IP algebraic curvature tensors are found, the corresponding
Riemannian metrics can be produced in a closed form, and depend on a few
constants. When $n =3$, the situation is completely different: the IP
algebraic curvature tensors can be easily classified \cite[Remark 1]{IP}:
they are either of constant curvature, or those whose Ricci tensor
$\rho$ has rank $1$. However, the class of Riemannian manifolds satisfying
the latter condition, $\rk \, \rho =1$, is quite large: it depends
on at least two arbitrary functions of one variable, and it seems doubtful
that the description of these manifolds can be obtained in some nice form.

As in dimension $3$ the Ricci tensor determines the curvature tensor, the
question of finding IP manifolds can be viewed as the question of finding a
Riemannian metric given its Ricci tensor. Even the existence of a solution $g$ for
the corresponding system of differential equations $\Ric(g) = \rho$
is a hard problem (see \cite[Ch.5]{Besse}
for examples of symmetric tensors which cannot be Ricci tensors of any Riemannian
metric). For nondegenerate Ricci tensors, the existence problem is solved in
affirmative by Deturck \cite{DT}. Recently, the existence of a Riemannian metric
$g$ with the given Ricci tensor $\rho$ was also proved for degenerate $\rho$ whose
kernel distribution has constant rank and is integrable (under some additional
assumptions on the first derivatives) \cite{DTG}.

Let $M^3$ be a Riemannian manifold whose Ricci tensor has constant rank one, with
$2f$ the nonzero principal Ricci curvature.

If $f = \mathrm{const}$ (and more generally, if the principal Ricci curvatures $\rho_i$
are constant and $\rho_1=\rho_2 \ne \rho_3$), the Riemannian manifold $M^3$, up to
isometry, depends on two functions of one variable, as was shown by Kowalski \cite{K,KS}
and Bueken \cite{Bu}. Despite of the fact that any such $M^3$ is curvature-homogeneous
(the curvature tensor at every point is the same), the majority of them are not homogeneous.
The only homogeneous $3$-manifolds with $\rk \,\rho=1$
are unimodular Lie groups with a specific left-invariant metric, whose explicit construction
is given by Milnor \cite[Ch.4]{M} (see also \cite{KN}). Let $\mathfrak g$ be a $3$-dimensional
Lie algebra with a basis $e_1, e_2, e_3$ and the Lie brackets defined by
$$
[e_1, e_2] = \la_3 e_3,\quad [e_2, e_3] = \la_1 e_1,\quad[e_3, e_1] = \la_2 e_2.
$$
Assuming $e_1, e_2, e_3$ orthonormal we get a left-invariant metric on the Lie group
$G$ of $\mathfrak g$. If $\la_1+\la_2 =\la_3, \; \la_1 \la_2 \ne 0$, the Ricci
tensor of this metric has rank $1$, with the nonzero principal Ricci curvature $2f = 2\la_1 \la_2$.
Depending on the signs of the $\la_i$'s, the underlying Lie group $G$ is $SU(2)$ (the sphere $S^3$,
but not with a constant curvature metric), $SL(2,\mathbb R)$, or $E(1, 1)$, the group of motions
of Minkowski plane. In the latter case, the metric has the form
\begin{equation}\label{eq:metricE11}
    ds^2 = dx^2 + e^{2ax}dy^2 + e^{-2ax}dz^2, \;a \ne 0,
\end{equation}
and is the only $3$-dimensional metric, which is generalized symmetric, but not symmetric
\cite[Ch.6]{Kg}.

If $f \ne \mathrm{const}$, only isolated examples are known.
We will consider here two particular cases: when the space $M^3$ is conformally flat, and
when the principal Ricci direction corresponding to $2f$ is a geodesic vector field.
This choice of the additional assumptions is motivated by the following facts. Firstly, all
the IP manifolds of dimension $n \ge 4$ (Example~\ref{ex:f(t)}) are conformally flat, which
is no longer true when $n = 3$. Secondly, when $f = \mathrm{const}$, the principal Ricci
direction corresponding to $2f$ is a geodesic vector field, which follows from the second
Bianchi identity
(see \eqref{eq:integrability} below).

\begin{proposition*}
Let $M^3$ be a Riemannian manifold whose Ricci tensor has rank one, with a nonconstant principal
Ricci curvature $2f$, and the corresponding principal Ricci direction $e_1$.
\begin{enumerate}
  \item If $M^3$ is conformally flat, then it is locally isometric to a manifold with metric \eqref{eq:f(t)}.
  \item If $e_1$ is a geodesic vector field, then $M^3$ is either conformally flat, or the metric form on
  $M^3$ is locally homothetic to
\begin{equation}\label{eq:metric}
  ds^2 = dx^2 + x^{1+a} dy^2 + x^{1-a} dz^2, \qquad \text{with} \, \; a \ne \pm 1.
\end{equation}
\end{enumerate}
\end{proposition*}


Before giving the proof, consider the general case. Let $M^3$ be a Riemannian manifold whose Ricci
tensor has constant rank one. Introduce a local orthonormal frame $e_1, e_2, e_3$ in such
a way that $e_1$ is the principal direction of the Ricci tensor corresponding to
$2f$, and $\Span(e_2, e_3) = \Ker \, \rho$. The only nonzero components of the curvature
tensor $R_{ijkl}=R(e_i,e_j,e_k,e_l)$, up to permutation of indices, are
$$
R_{1212}= R_{1313}= f, \quad R_{2323}=-f.
$$
Let $\om^i$ be the $1$-forms dual to $e_i$, and let $\psi_i^j, \, \Om_i^j$
be the connection and the curvature forms, respectively:
$$
\psi_j^i = \G_{jk}^i \om^k, \;\G_{jk}^i = \<\n_k e_j, e_i\>, \;\G_{ik}^j = - \G_{jk}^i ,
\;\psi^i_j =- \psi^j_i, \; \Om^i_j = - \Om^j_i = \tfrac 12 R_{ijkl} \om^k \wedge \om^l,
$$
so that
\begin{equation}\label{eq:curvforms}
\Om^1_2 = f \om^1 \wedge \om^2, \qquad \Om^1_3 = f \om^1 \wedge \om^3,
\qquad \Om^2_3 = - f \om^2 \wedge \om^3.
\end{equation}
We have the structure equations
\begin{equation}\label{eq:structure}
d \om^i = - \psi^i_j \wedge \om^j, \qquad d \psi^i_j =
- \psi^i_k \wedge \psi^k_j + \Om^i_j,
\end{equation}
whose integrability condition is the second Bianchi identity
\begin{equation}\label{eq:Bianchi}
d \Om^i_j = \Om^i_k \wedge \psi^k_j - \Om^k_j \wedge \psi^i_k.
\end{equation}
Substituting \eqref{eq:curvforms} to \eqref{eq:Bianchi} we find:
\begin{equation}\label{eq:integrability}
  \frac{e_1(f)}{2f}= \G_{22}^1  + \G_{33}^1, \quad \frac{e_2(f)}{2f} = \G_{11}^2, \quad \frac{e_3(f)}{2f}= \G_{11}^3,
\end{equation}
or, equivalently, $df/(2f)= (\G_{22}^1  + \G_{33}^1) \om^1 + \G_{11}^2 \om^2 + \G_{11}^3 \om^3$.
Another equivalent form of the second Bianchi identity is
\begin{equation}\label{eq:integrability+}
  d(\sqrt{|f|} \om^1) = \sqrt{|f|} (\G_{23}^1  - \G_{32}^1) \om^2\wedge \om^3, \qquad d(\sqrt{|f|}\,\om^2\wedge \om^3) = 0.
\end{equation}

Studying the system (\ref{eq:curvforms}, \ref{eq:structure}, \ref{eq:integrability+}) further,
it might be interesting to know, for example, whether the solution set (say, in the analytic
case) depends on functions of two variables.

\begin{proof}[Proof of the Proposition]
1. A manifold $M^3$ is conformally flat, if its Schouten-Weyl tensor vanishes, that is,
if the tensor $T(X, Y, Z) = (\n_X \rho)(Y,Z) -\tfrac14 \<Y, Z\> X(s)$ is
symmetric with respect to $X, Y$, where $s$ is the scalar curvature. In
our case, $\rho = 2f \om^1 \otimes \om^1, \; s = 2f$. A direct calculation
shows that $M^3$ is conformally flat, if and only if
\begin{equation}\label{eq:confflat}
    e_1(f) = 4 \G^1_{22}f = 4 \G^1_{33}f, \quad e_2(f) = e_3(f) = 0, \quad \G_{11}^2=\G_{11}^3=0
    , \quad \G_{23}^1=\G_{32}^1=0,
\end{equation}
and then the first Bianchi identity \eqref{eq:integrability} is automatically satisfied.

From this point on, the proof goes word-by-word as in the four-dimensional case, starting
from equation~(3.21) on page~$279$ of \cite{IP}, up to changing the notation.

\vskip .3cm
\noindent2. As the field $e_1$ is geodesic, $d\om^1 = 0$, and we can choose a (coordinate)
function $x$ on a neighbourhood $U \subset M^3$ in such a way that $\om^1 = dx$.

The proof goes in seven steps:

\noindent\textbf{Step 1.} The distribution $\Ker \, \rho = \Span(e_2, e_3)$ is integrable and $f = f(x)$.

As $e_1$ is geodesic, $\G_{11}^2 = \G_{11}^3 = 0$, so by \eqref{eq:integrability}, the
fields $e_2, e_3$ are tangent to the level sets of $f$. Again, from \eqref{eq:integrability},
$df$ is a scalar multiple of $\om^1=dx$, so $f$ is a function of $x$. As $\Span(e_2, e_3)$ is
integrable, $\G^1_{23} = \G^1_{32}$.

\noindent\textbf{Step 2.} The fields $e_2, e_3$ can be chosen in such a way that $\G_{31}^2 = 0$.

Replacing $e_2, e_3$ by $\tilde e_2 = \cos \a \, e_2 + \sin \a \,e_3, \; \tilde e_3=-\sin \a \,e_2 + \cos \a \,e_3$,
respectively, with some function $\a$, we find that
$\tilde \G_{31}^2 =\<\n_1 \tilde e_3, \tilde e_2\> = \G_{31}^2 - e_1(\a)$. Choosing $\a$ in such a way that
$e_1(\a) = \G_{31}^2$ we get what required.

Let $H$ be a symmetric $2\times 2$ matrix, with entries $h_{ij}= \G^1_{ij},\; i,j = 2,3$ (the second fundamental form
of the foliation $f = \mathrm{const}$). As it follows from \eqref{eq:integrability},
\begin{equation}\label{eq:TrH}
  \Tr \, H = f^\prime/(2f).
\end{equation}

\noindent\textbf{Step 3.} The matrix $H$ satisfies the differential equation
\begin{equation}\label{eq:dH}
e_1(H) = H^2 + fI_2.
\end{equation}

From the structure equations \eqref{eq:structure}, the form of the
curvature tensor \eqref{eq:curvforms}, and the fact that
$\psi^1_i = h_{ij} \om^j, \; i = 2, 3$, we obtain
$$
(dH - \om^1 (H^2 + fI_2) - \psi^2_3\,[H, J])\wedge\left(\begin{matrix} \om^2\\ \om^3 \end{matrix}\right) =0,
$$
where $J = \left(\begin{smallmatrix} 0&1\\-1&0\end{smallmatrix}\right)$.
Extracting the $\om^1\wedge\om^i$ components from the both rows, and using the fact that
$\G_{31}^2 = 0$ (Step 2), we get \eqref{eq:dH}.

\noindent\textbf{Step 4.} Let $f(x)=\varepsilon \phi^{-4}(x), \; \varepsilon = \pm 1$. Then
\begin{equation}\label{eq:hij}
h_{22} =\phi^{-2} v - \phi^{-1} \phi', \quad h_{33} =-\phi^{-2} v - \phi^{-1} \phi', \quad h_{23} =\phi^{-2} u,
\end{equation}
where the functions $u, v$ satisfy $e_1(u) = e_1(v) =0, \; u^2 + v^2 + \varepsilon = c_0 = \mathrm{const}$, and
\begin{equation}\label{eq:phi}
\phi^2 =A x^2 + Bx +C, \qquad \text{with constants $A, B, C$ satisfying}\; B^2-4AC = 4c_0.
\end{equation}

This can be obtained directly by solving the system of ODE's (\ref{eq:TrH}, \ref{eq:dH}).

\noindent\textbf{Step 5.} The fields $e_2, e_3$ can be chosen in such a way that $u = 0$, and both
distributions $e_2^\perp,\, e_3^\perp$ are integrable.

Replacing $e_2, e_3$ by $\tilde e_2 = \cos \b \, e_2 + \sin \b \,e_3, \; \tilde e_3=-\sin \b \,e_2 + \cos \b \,e_3$,
respectively, with some function $\b$ such that $e_1(\b) = 0$ (not to violate the condition of Step 2), we find:
$\tilde \G_{23}^1 =\<\n_{\tilde 3} \tilde e_2, e_1\> = \cos 2\b\, \G_{23}^1 + \tfrac12 \sin 2\b \,(\G_{33}^1 - \G_{22}^1)=
\cos 2\b\, h_{23} + \tfrac12 \sin 2\b \,(h_{33} - h_{22})= \phi^{-2} (\cos 2\b\, u - \sin 2\b \, v)$ by \eqref{eq:hij}.
Choosing $\b$ accordingly, we obtain $\tilde \G_{23}^1 = 0$.

Omitting the tildes, we get $\G_{23}^1 = \G_{32}^1 = h_{23}=u=0$. Then
$d\om^2 \wedge \om^2 = (\G^2_{31}-\G^2_{13}) \, \om^1\wedge\om^2\wedge\om^3 = 0$ (from the above and Step 2),
and similarly $d\om^3 \wedge \om^3 = 0$. Note also that $v^2 = c_0 -\varepsilon = \mathrm{const}$.

\noindent\textbf{Step 6.} $v \, \psi^2_3 = 0$.

We already know that
$\psi_2^1 = (\phi^{-2} v - \phi^{-1} \phi')\,\om^2,\; \psi_3^1 = (-\phi^{-2} v - \phi^{-1} \phi')\,\om^3,\;
\psi_3^2 = \G^2_{32} \om^2 + \G^2_{33} \om^3$. Substituting this to the structure equation
$d\psi^1_2 = -\psi^1_3 \wedge \psi^3_2 + f \om^1 \wedge \om^2$ and extracting the $\om^2 \wedge \om^3$
term we get $v \G^2_{32} = 0$. Similarly, $v \G^2_{33} = 0$.

\noindent\textbf{Step 7.} Assuming $v = 0$, we go back to the conformally flat case. Indeed,
equation \eqref{eq:hij} implies $h_{22} = h _{33} =- \phi^{-1} \phi' = e_1(f)/(4f)$, also
from Step 1 we know that $e_i(f) = 0, \G_{11}^i=0, \; i =2, 3$, and
from Step 5, $\G_{23}^1 = \G_{32}^1 =0$. Then \eqref{eq:confflat} follows.

Let us take $v \ne 0$. By Step 6, $\psi^2_3 = 0$. It follows that
$d\om^2 = -(\phi^{-2} v - \phi^{-1} \phi')\,\om^1 \wedge \om^2 =
-(\phi^{-2} v - \phi^{-1} \phi')\, dx \wedge \om^2$, so we can find functions $\mu_2 = \mu_2(x)$
and $y$ such that $\om^2 = \mu_2(x) dy$. Similarly, for some functions $\mu_3(x)$ and
$z, \; \om^3 = \mu_3(x) dz$. Then
$$
ds^2 = (\om^1)^2+(\om^2)^2+(\om^3)^2 = dx^2 + \mu_2(x)^2 dy^2 + \mu_3(x)^2 dz^2.
$$
Calculating the Ricci tensor
(with \textsc{Maple}), and equating $\rho_{22}$ and $\rho_{33}$ to zero, we get
$\mu_2'' \mu_3 + \mu_2' \mu_3' = \mu_3'' \mu_2 + \mu_2' \mu_3' = 0$, so, up to homothecy and translation,
$$
\mu_2 = x^{(1+a)/2},\; \mu_3 = x^{(1-a)/2},\; \qquad\text{or}\qquad \mu_2 = e^{ax}, \; \mu_3 = e^{-ax}.
$$
In the second case, we get the metric form \eqref{eq:metricE11}, with
$2f = \rho_{11} = -2a^2 = \mathrm{const}$, which contradicts the assumption.

The first case gives the required metric \eqref{eq:metric}. We have
$2f = \rho_{11} = \frac12 (1-a^2) x^{-2}$,
and the metric \eqref{eq:metric} is not isometric to any of \eqref{eq:f(t)}, unless $a=0$
(for instance, because the surfaces $f = \mathrm{const}$ are not totally umbilical).
\end{proof}

We finish with yet another example of a metric whose Ricci tensor has rank one, and the scalar
curvature is nonconstant:
$$
ds^2=e^y dx^2 + \frac{e^y}{y}dy^2+ y dz^2.
$$
The nonzero principal Ricci curvature of this metric is $-\frac12 e^{-y}$, with the
corresponding principal Ricci direction $e^{-y/2}\partial/\partial x$ (which is
not geodesic).



\begin{thebibliography}{99}

\bibitem{Besse}
A.Besse, Einstein manifolds,
Ergebnisse der Mathematik und ihrer Grenzgebiete (3),
\textbf{10}. Springer-Verlag, Berlin, 1987.

\bibitem{Bu}
P.Bueken,
Three-dimensional Riemannian manifolds with constant principal Ricci
curvatures $\rho\sb 1=\rho\sb 2\not=\rho\sb 3$,
J. Math. Phys.
\textbf{37} (1996), 4062--4075.


\bibitem{DTG}
D.DeTurck, H.Goldschmidt,
Metrics with prescribed Ricci curvature of constant rank. I. The
integrable case,
Adv. Math.
\textbf{145}, (1999), 1--97.

\bibitem{DT}
D.DeTurck,
Existence of metrics with prescribed Ricci curvature: local theory,
Invent. Math.
\textbf{65}, (1981/82), 179--207.

\bibitem{GLS}
P.Gilkey, J.Leahy, H.Sadofsky,
Riemannian manifolds whose skew-symmetric curvature
operator has constant eigenvalues,
Indiana Univ. Math. J.
\textbf{48}, (1999), 615--634.

\bibitem{G2}
P.Gilkey,
Riemannian manifolds whose skew-symmetric curvature
operator has constant eigenvalues II,
in: Differntial Geometry and Applications, Masaryk Univ, Brno, 1999, 73 --87.

\bibitem{GS}
P.Gilkey, U.Semmelmann,
Spinors, self-duality, and IP algebraic curvature tensors,
ESI preprint 616, (1998).

\bibitem{Gb}
P.Gilkey,
Geometric properties of natural operators defined by the Riemann curvature tensor,
World Scientific Publishing Co., Inc., River Edge, NJ, 2001.

\bibitem{GI1}
P.Gilkey, R.Ivanova,
The geometry of the skew-symmetric curvature operator in the complex setting,
Global differential geometry: the mathematical legacy of Alfred Gray (Bilbao, 2000),
Contemp. Math., \textbf{288}, 325--333.

\bibitem{GI2}
P.Gilkey, R.Ivanova,
Complex IP pseudo-Riemannian algebraic curvature tensors,
PDEs, submanifolds and affine differential geometry (Warsaw, 2000),
Banach Center Publ., \textbf{57}, Polish Acad. Sci., Warsaw, (2002),  195--202.

\bibitem{GZ}
P.Gilkey, T.Zhang,
Algebraic curvature tensors for indefinite metrics whose skew-symmetric curvature operator has
constant Jordan normal form,
Houston J. Math. \textbf{28} (2002), 311--328.

\bibitem{IP}
S. Ivanov, I.Petrova, Riemannian manifold in which the skew-symmetric curvature
operator has pointwise constant eigenvalues,
Geom. Dedicata
\textbf{70}, (1998), 269--282.

\bibitem{Kg}
O. Kowalski,
Generalized symmetric spaces,
Lecture Notes in Mathematics,
\textbf{805}, Springer-Verlag, Berlin-New York, 1980.


\bibitem{K}
O.Kowalski,
A classification of Riemannian $3$-manifolds with constant
principal Ricci curvatures $\rho\sb 1=\rho\sb 2\not=\rho\sb 3$,
Nagoya Math. J.
\textbf{132} (1993), 1--36.

\bibitem{KS}
O.Kowalski, M.Sekizawa,
Local isometry classes of Riemannian $3$-manifolds
with constant Ricci eigenvalues $\rho\sb 1=\rho\sb 2\neq \rho\sb 3>0$,
Arch. Math. (Brno) \textbf{32} (1996), 137--145.

\bibitem{KN}
O.Kowalski, S. Nik\v cevi\'c,
On Ricci eigenvalues of locally homogeneous Riemannian $3$-manifolds,
Geom. Dedicata
\textbf{62} (1996), 65--72.

\bibitem{M}
J.Milnor,
Curvatures of left invariant metrics on Lie groups,
Advances in Math.
\textbf{21} (1976), 293--329.

\bibitem{Nag}
M.Nagata, A remark on the unique factorization theorem,
J. Math. Soc. Japan,
\textbf{9}, (1957), 143--145.

\end{thebibliography}
\end{document}